# Environmental Impacts of MSW Collection Route Optimization Using Gis: A Case Study Of 10th Of Ramadan City, Egypt

Amr Shafik, Moustafa Elkhedr, Dalia Said, Ahmed Hassan

**Abstract**

For municipal solid waste (MSW) collection purposes, optimal routes for collection need to be determined for trucks to improve the impact on both the environment and the collection cost. This improvement leads to benefits such as decreasing truck emissions, and pollution and improving the collection efficiency of the process. In this research, a case study is developed using the ArcMap Network Analyst tool to plan the MSW collection and transportation scheme from the city of the 10$^{th}$ of Ramadan to a new proposed landfill location. Several collection scenarios are analyzed and evaluated in terms of total travel time and distance to find the best operation variables such as the number of collection trucks and their routes. Results showed that the system improved the total travel time of waste collection by (26.5%). However, due to the increased number of trucks, the traveled distance increased by (90%) in the proposed scenario.

## 1. Introduction

Municipal solid waste (MSW) forms major environmental, social, and economic impacts on urban communities, especially with the continuous increase in population. The composition and generation rates of MSW are affected by the geographic location of the area being studied, population, income level, and other socio-economic conditions (Akinci et al. 2012; Al-jarallah and Aleisa 2014; Chandrappa and Das 2012; Magrinho et al. 2006; Tınmaz 2006). MSW produced in higher-income countries is usually of more value than in lower-income countries. Higher-income regions usually generate recyclable materials such as cardboard, glass, and plastic with an annual growth rate of 3.2 – 4.5%. MSW collected in low-income regions is mostly degradable material with an annual growth rate of 2-3% (Lu et al. 2015; Shekdar 2009).

The MSW management process includes waste collection, transfer, treatment, and disposal. Cost-effective solutions for each component are vital for a successful integrated system. The cost component of MSW collection and transport to treatment facilities or transfer stations alone may reach 60-70% of the system cost and is affected by the quantity and type of waste, as well as the method of transport. MSW management in developing countries is usually based on experience and lacks efficiency in optimizing resources. Since any improvement in this component would have an overall positive impact on the system, thorough planning for efficient and effective use of resources and costs is essential to cut unnecessary costs.

Extensive research has been conducted to optimize the MSW management process. Improvement strategies of MSW collection include optimization of MSW truck routing to decrease fuel consumption and traveled distance, and consequently less greenhouse gas emissions and increased life span of used equipment. Additional environmental benefits of optimized vehicle routing include decreased contamination of environmental elements at waste disposal sites. In recent years, due to the significant increase in computational power, optimization techniques are now possible and have been used as a decision support tool for problem-solving in many disciplines including waste management. Optimization methods include techniques that are based on Geographic Information Systems (GIS). Although optimization techniques have been used widely in developed regions, there is a shortage in their application in developing countries. The objective of this paper is the optimization of MSW collection and transfer costs using GIS as an optimization tool by tuning the optimal number of trucks and finding their minimum-cost routes that minimize the travel time required to perform the required task of waste collection from given locations.

Egypt is a middle-income country and is one of the largest waste-generating countries with an average current waste generation rate estimated at 0.58 kg/capita/day in urban areas with a growth rate of 2% and an approximate daily generation rate of 57,000 tons. The Greater Cairo Region alone generates approximately 34% of this waste. The case study is selected to be 10th of Ramadan City, a new urban community in Egypt East of Cairo, and generates approximately 565 tons of waste per day. MSW in the city is currently based on private waste collectors.

The current strategy of the Ministry of Environment (MOE) and the New Urban Communities Authority in Egypt is the development of waste management masterplans in new urban communities to prevent the negative impacts of unplanned waste removal. The MOE has put in place waste management masterplans in new urban communities East, West, and South of the Greater Cairo Region with a total of 8 cities East of GCR, 7 cities west of GCR, and 11 cities south of GCR. Part of the master plan was to optimize MSW collection and transfer costs using GIS.

Local municipalities are responsible for waste collection from households. The municipalities, in turn, contract private waste collectors to collect the waste from households (door–to–door) and take a monthly fee from each household. Private waste collectors use light open trucks to transport waste to transfer stations, and then large transfer trucks to transport the waste to treatment facilities or landfills. Factors influencing the cost include MSW quantity, number and capacity of waste collection trucks, and the shortest cycle path to reduce fuel costs. Currently, the routing system is determined by the experience of the waste collectors and not by an optimization technique.

The objective of this study is to use a GIS-based optimization method to minimize the collection and transfer cost for a solid waste collection system for the study area with a proposed treatment and disposal site; and therefore, reduce the travel time and fuel consumption on the network to reach to the minimum greenhouse gas emissions. The used optimization technique used in the GIS platform is based on Dijkstra's Algorithm which is a minimum-cost path finding heuristic. In doing so, the cost components were well defined, and the variables that could be minimized were determined. The operational parameters, including the number of collection trucks, the size of labor in each truck, and truck routes, were applied and evaluated to find the optimal scenario for the operation.

Cost components in the collection and transfer phase include fixed costs and variable costs of vehicles and costs related to labor, equipment, and overhead costs (Elkhedr 2016). Since the fixed costs and expenditures related to labor and equipment are not variable, they were not included in the optimization technique. Variables that could be optimized are the variable costs that depend on the distances traveled by collection vehicles including fuel costs, repair and maintenance costs, tire costs, as well as the fleet size, number of trips taken per day to the transfer station, truck capacity, and average number of apartments served per trip.

## 2. Background

### 2.1 MSW Collection Resources and Techniques

As the collection of Municipal Solid Waste (MSW) is the main step in the waste management system, governments have developed various collection strategies that depend on several factors, such as resident types, community characteristics, available resources, and types of the collected waste itself. Each strategy implements available logistic resources towards achieving the required objectives with minimum resources, such as waste collection bins and containers, pickup locations, landfill site locations, and collection and transportation equipment (Whitman 2005). Effective systems are also manifested in reducing the operation cost of these resources, as well as their negative environmental impacts. There are two types of waste collection systems presented in the literature; hauled container systems (HCS), and stationary container systems (SCS) (Theisen and Vsigil, 1993). In HCS, waste storage bins are moved from their original positions to the disposal sites to be emptied and returned afterward. While in SCS, the waste storage bins are emptied at their locations. Either of the two types can be effectively used depending on the community and resource characteristics as mentioned before.

The "Network Analyst" tool which is available in "Esri ArcMap" software on the GIS platform was used for this application, this tool is based on Dijkstra's algorithm which finds the shortest path between any set of nodes in a graph. In this application, the travel time would be the objective function to minimize. Although several studies have discussed this application, as presented in Table 1, it has not been applied to the unique case study in Egyptian urban communities, and research in this area using GIS optimization techniques is somewhat new. The traditional vehicle routing problem formulation which minimizes the total traveling cost, which is considered the travel time, in this case, is commonly found in the literature. The formulation found in (Kulkarni & Bhave, 1985) is the one used in this analysis.

**Table 1: Previous Work on MSW Collection and Transportation Optimization Using GIS Techniques**

| Application Country | Contribution | Reference |
|---|---|---|
| Port Said, Egypt | Optimization of municipal solid waste management using mixed integer programming | (Badran 2006) |
| Abeokuta, Nigeria | Identification of disposal sites using GIS and RS | (Ufoegbune 2016) |
| Santiago, Chile | Establish feasible collection routes - determine vehicle fleet size - present comparative cost and sensitivity analysis of the results | (Arribas et al. 2010) |
| Laxmi Nagar, Nagpur, India | design and develop an appropriate collection plan by GIS, optimize existing collection cost | (Bhambulkar & Khedikar 2011) |
| India | Optimize MSW collection routes length | (Das & Bhattacharyya 2015) |
| Tanzania | Examining how GIS can increase efficiency and info of solid waste collection systems | (Kyessi & Mwakalinga 2016) |
| Indonesia | Assessment of waste collection routes in ArcGIS, LOS Assessment | (Hasyim et al. 2018) |
| Jordan | Optimization of working hours - routes - transfer stations | (Oelgemöller & Nelles 2017) |
| Kampala, Uganda | Optimize travel distance, trips, and collection time - Identify the optimum location of the proposed landfill site | (Kinobe et al. 2015) |
| Turkey | optimize solid waste collection/hauling processes, as the minimum cost was aimed - the optimized routes were compared with the present routes. | (Apaydin 2014) |
| Sweden | develop a model to calculate time and energy consumption | (Sonesson 2000) |
| Portugal | Develop a system for route optimization and collection scheduling | (Zsigraiova, et al. 2013) |
| UAE | Route optimization and comparison with the actual situation | (Maraqa, et al. 2018) |
| Ghana | Determine the effect of route optimization on travel distance, travel time, and fuel consumption of municipal solid waste (MSW) collection trucks | (Sulemana, et al. 2019) |
| Serbia | An optimization methodology for collection and transportation of municipal solid waste in large urban centers - calculation of fuel consumption, and optimization of waste collection and transportation using GIS | (Ristić et al. 2015) |

## 2.2 Selection of Collection Stop Points and Disposal Sites

The selection of the collection stop points' locations is a crucial question as it has a significant effect on the efficiency of the collection process, accessibility of the collection equipment, and the collection operation time and cost (Lu et al. 2015). Careful selection of the optimal locations is addressed in literature by several approaches. For example, a Genetic Algorithm-based Multi-Objective Integer Programming was used to allocate collection points in the collection network to optimize collection time and cost (Chang and Wei 2000). Moreover, similar objectives were achieved using Binary programming (Kao and Lin 2002).

Regarding the selection of disposal and landfill sites, a GIS analysis-based Integer Programming approach (IP) was used to determine the optimal locations of transfer stations in a large metropolitan area (Chang and Lin, 1997). Another research solved the problem of choosing disposal facilities heuristically to achieve an optimal operation system (Benjamin and Beasley 2013). In addition, Dey et al. (2019) used a GIS-based multi-criteria analysis technique to evaluate the suitability of three existing landfill sites and to evaluate several relevant parameters that are important for landfill sites.

## 2.3   Vehicle Routing Problem (VRP)

The waste collection and transportation process has several resources and equipment that need to be optimized and managed appropriately to achieve the intended objectives. These resources include collection facilities, such as collection locations, landfills, and treatment sites. Besides, these include collection equipment like collection and transportation trucks and labor resources. All these elements require an optimization system to achieve the optimal balance between demand and supply while taking into consideration other constraints like a limited time window, limited operation and maintenance cost, environmental impacts of operation processes such as trucks emissions, and community contamination due to the presence of uncollected waste for some time.

Various approaches were presented in the literature to solve the vehicle routing problem of waste collection and transportation. Das and Bhattacharyya (2015) formulated the VRP into a mixed-integer program. The optimal waste collection and transportation system was found by implementing the created integer program together with a proposed heuristic solution. Another research implemented the Clark and Wright algorithm to find optimal vehicle routes (Sarmah, Yadav, and Rathore 2019). Moreover, some studies used the Particle Swarm Optimization algorithm (PSO) to find optimal routes (Hannan et al. 2018; Son 2014). In addition to these approaches, some researchers used Artificial Neural Networks to assess waste collection systems (Vu et al. 2019). Others used an agent-based model together with GIS (Hua et al. 2017). Finally, a comprehensive review of mathematical programming approaches to perform the task of optimizing waste collection trucks was presented in Sulemana et al. (2018).

## 2.4 Using GIS to Solve Vehicle Routing Problem (VRP)

Geographic Information System (GIS) is commonly used in the process of route planning for MSW collection and transportation trucks. Some researchers used GIS to improve the current routes of waste collection. Maraqa et al. (2018) developed a GIS-based model for Um Gafa in the City of Al Ain, United Arab Emirates, to calculate the fuel consumption rate for the process of MSW collection. The model was used to reduce the cost of the waste collection process, as well as the greenhouse emissions by optimizing the fuel consumption during MSW collection. The methodology consisted of three stages: developing a GIS-based Model, calculating fuel consumption rates for different scenarios, and investigating the adequacy of waste collection bins in the study area. Results showed that when using the optimization approach, fuel consumption decreased by 14.3%. In addition, greenhouse gases such as $CO_2$, $CO$, $HC$, and $NO_x$ were significantly decreased.

Moreover, Sulemana et al. (2019) performed route optimization using GIS to determine its effect on the travel distance, travel time, and fuel consumption of municipal solid waste (MSW) collection trucks. In this study, the network analyst extension in ArcMap was used to model and analyze the transportation network for the chosen study area in Ghana. Results showed that the travel distance was reduced by 4.79% when using the optimal routes. Also, the reduction in travel time and fuel consumption was 14.21% and 10.81% respectively, significantly reducing the operation cost of the MSW collection process.

Moreover, the GIS vehicle routing problem was widely used in many case studies such as the work of Chalkias and Lasaridi (2018), which improved the system of waste collection and transportation in the Municipality of Nikea, Athens using GIS technology. Ufoegbune(2016) applied the use of GIS and remote sensing to select waste disposal locations and routing analysis in Abeokuta, Nigeria.

## 3. Methodology

Trucks that collect MSW generated from households have different possible routes to collect from a series of households. Based on the chosen route, the time taken and fuel consumption would vary. The proposed methodology of this study included four steps as shown in Figure 1. Step 1 was collecting spatial data on the characteristics of the chosen study area and the existing waste collection method. Step 2 was setting up the GIS platform with data related to the road network and its characteristics, vehicles used for collection, and information on land use and population densities to determine MSW generation rates and quantities. Step 3 was the application of the GIS Network Analyst tool to find the optimum stop points, optimum number of trucks, and optimum vehicle routing. Finally, Step 4 was devoted to quantifying the benefits of using this methodology in terms of cost and time, and therefore translated to environmental benefits.

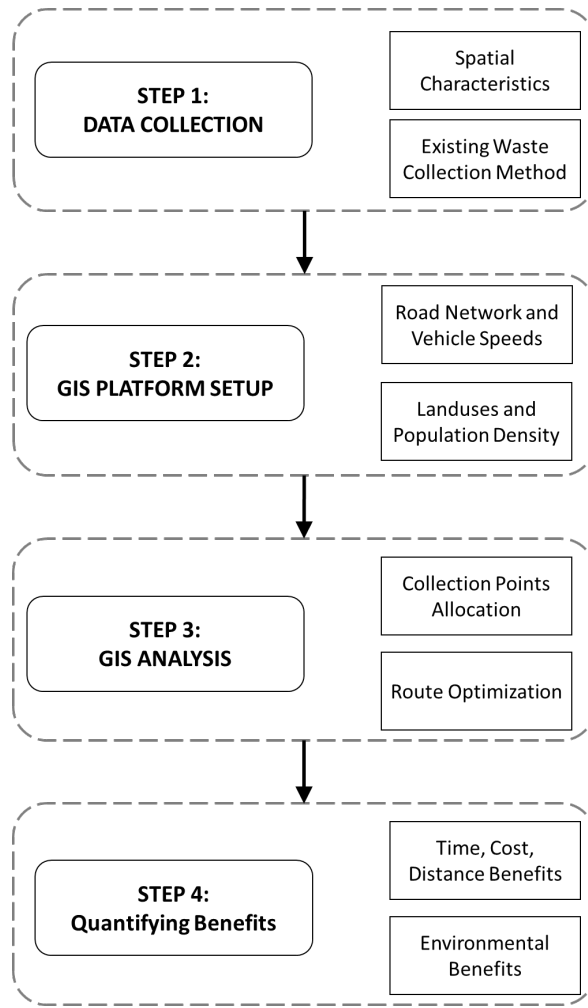

**Figure 1: Research Methodology**

## 3.1 Data Collection

The methodology of this study was applied to the 10$^{th}$ of Ramadan City, east of Cairo, Egypt. This study focuses on MSW collected from residential areas (61% of land use) as this waste is collected separately. The trucks collecting the MSW will be designed to stop at designated collection points to collect waste directly from households, then transport the MSW to a proposed 840,000 m² waste treatment and disposal facility located south of the city. The methodology involved designating specific stop points and rescheduling collection times to minimize route distances and collection time and costs.

## 3.2 Study Area

10th of Ramadan City is considered one of the new cities over 29 km² and consists of residential, industrial, recreational, and institutional land uses as shown in Figure 1. The population is approximately 850,000 capita with 197,000 households. Family size ranges from 4 - 5 with an average of 4.32. Most households reside in multi-story buildings including from 8 – 10 apartments per building in mixed residential and commercial land use.

10th of Ramadan City currently generates 565 tons/day with a generation rate of 2.49 kg/unit/day. The current waste collection system is operating using 16 trucks with an 18-ton capacity. MSW is collected from dumpsters with the stationary container systems (SCS) in the city and transported to an 840,000 m² waste treatment and disposal facility located south of it. In the proposed system, the crew size is 3 laborers/truck, one driver, and two collectors that collect from households and dump in the truck. Trucks may do more than one trip per day to the city and the MSW collection interval is every day, and the fleet size will be calculated from the optimization technique.

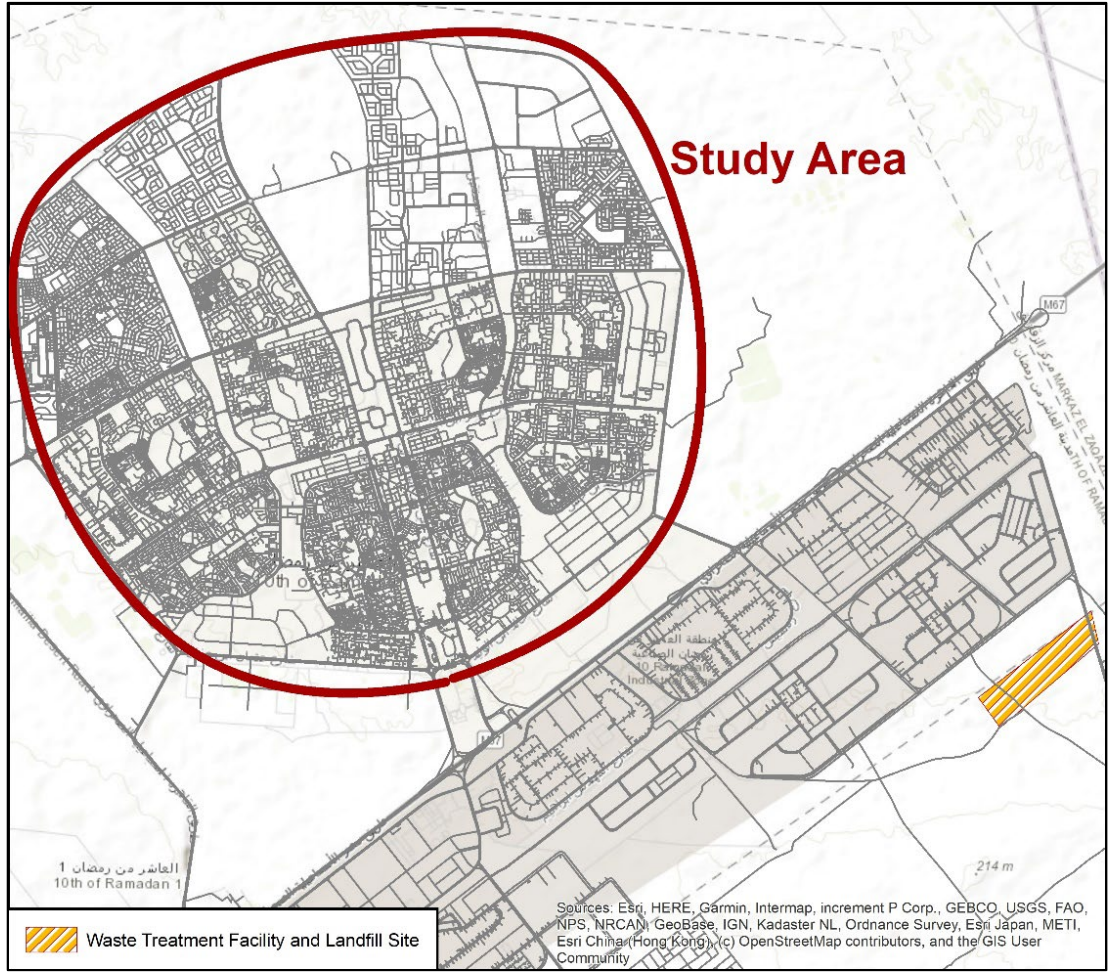

**Figure 2: Waste Treatment and Disposal Facility**

### 3.3 GIS Platform Setup

The optimization problem is analyzed using the spatial analysis capability of GIS, where the ArcMap software is set as a platform for all collected data. The database saved in the GIS database included attributes such as satellite image, land uses, population, the road network, its attributes, and existing stop points.

The platform was set up with data related to the road network and MSW data. The road network information and existing traffic volumes were collected from site visits, traffic surveys, and online GIS datasets. Geometric road characteristics such as road lengths, capacities, number of lanes, and speed limits were compiled into the GIS database. In addition, existing traffic volumes and classified counts were collected through traffic surveys at selected seven locations on the network on a typical working day to determine peak hour volumes and vehicle composition. This helps in determining optimum collection times to avoid any additional effect on the current traffic system.

The MSW required data included waste generation rate, truck capacities, service areas, and service schedule. In addition, population data, land uses, collection routes, truck speeds, fuel consumption, and emissions data were compiled into the GIS platform.

## 3.4   GIS Analysis

### 3.4.1   Optimized Stop Points

The optimal stop points locations, where the collection trucks will stop and wait for the crew to collect the waste from households, are set inside the study area based on several factors/inputs including:

- The collection method is performed door-to-door.
- Each stop point duration is 30 minutes.
- The collection workers can cover a service area of a radius of 300 m around the stop point. The 300m-radius circle contains on average approximately 26 residential buildings. Each building is composed of 4 floors, with 2 dwelling units per floor.

Given these inputs, trucks collect approximately 500 kg of MSW per stop point. Stop points were optimized and determined using the Network Analyst module in ArcMap software by ensuring that the collection locations achieve the coverage criteria mentioned before. The results of the optimization process of stop points showed that the required number is 500 stop points. This is opposed to the current 381 dumpster points used in the city.

Additional input data defined in the model were collection schedule, time windows, truck data, collection locations data, and amount of collected solid waste. Therefore, the whole study area was covered with the appropriate number of stop points that is sufficient for the quantity of daily generated waste and served the required area as shown in Figure 3.

### 3.4.2   Optimized Number of Trucks for MSW Collection

Based on the required number of stops and their locations, the number of collection trucks for the study area was defined using Network Analyst Tool. In addition, other constraints control the number of trucks, such as the truck's capacity, working time window, and the waiting time of each truck at stop points. The truck's capacity for the proposed collection scheme is 4 tons. The waiting time for trucks (collection time) is 30 minutes at each stop point. Finally, the average speed of trucks on the road network is 40 km/hr.

The analysis is performed considering the given constraints by the Network Analyst. The tool in ArcMap ensures that these trucks traverse the given set of stop points with minimum travel time and travel distance. The network analyst tool analyses the number of trucks in the system and evaluates possible scenarios to find the best collection scenario for the city and the optimal number of trucks.

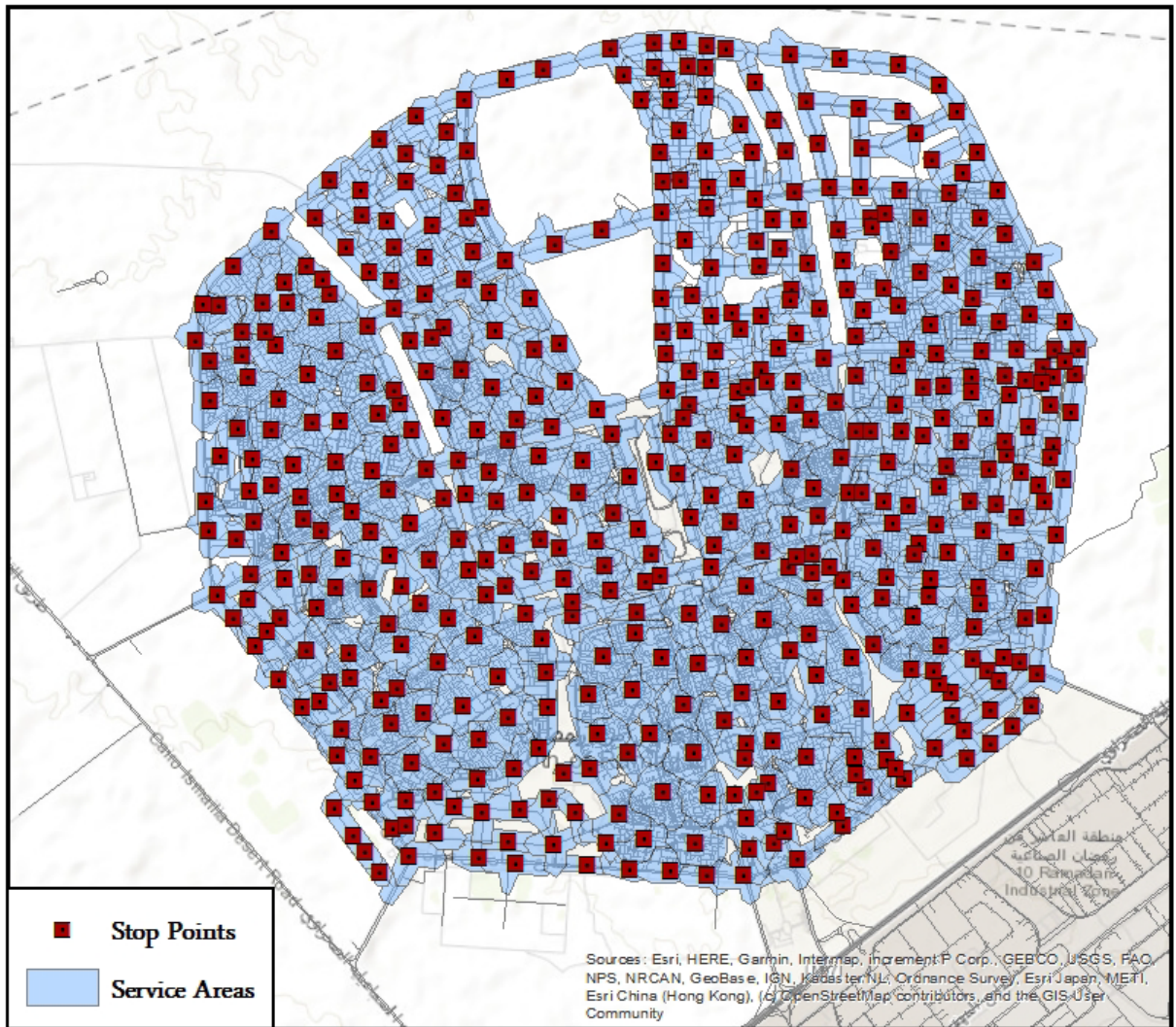

**Figure 3: Proposed Stop points and Service Areas around each Collection Point (Radius<300m)**

### 3.4.3 Optimized MSW Trucks Routing

Optimal routes of each specified vehicle were generated in the Network analyst tool using Dijkstra's Algorithm which finds the optimal path to cover the required locations with minimum cost. In Dijkstra's algorithm, a single-source shortest path is found between any pair of nodes on a weighted directed graph (Cormen et al. 2009). The algorithm is applied in this study to find the optimal routes for waste collection trucks represented in the shortest path between each successive pair of stop points. These routes are also constrained by trucks' capacities and working time windows. Each collection truck serves as a single problem in the entire system of trucks, stop points and the road network. The optimization problem finds the routes for each truck (per cycle) to cover all the stop points with minimum overall travel time (or travel distance). The modeled geometric characteristics of the road turns and intersections represented actual field characteristics.

Figure 4 shows the generated trucks' routes by the Network Analyst Tool. During the assignment process, the trucks traverse the multiple stop points until it reaches their full capacity. Then, they return to the treatment and disposal facility location to unload and start over again in another route for waste collection.

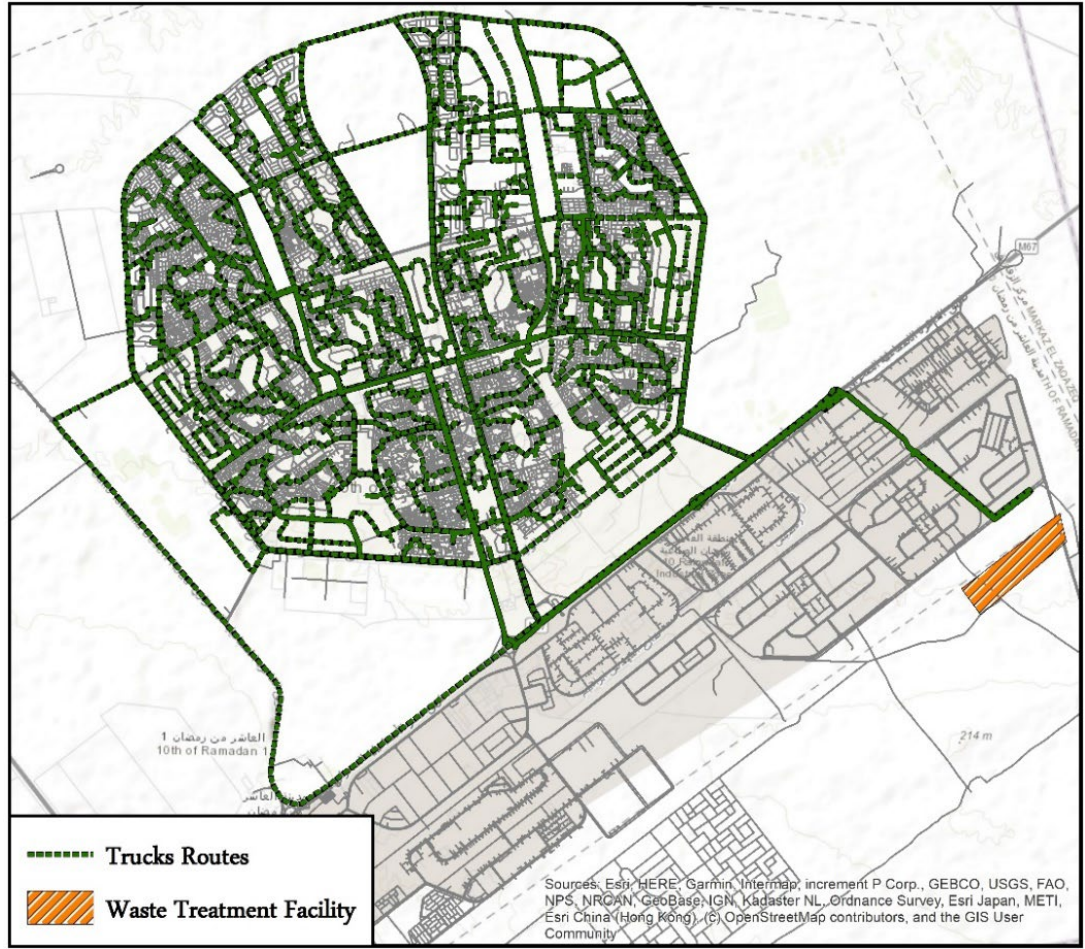

**Figure 4: Generated Optimal Trucks Routes**

## 4. Results and Discussion

### 4.1 Quantification of Cost and Time Impact

Table 2 shows a comparison between the different elements of the existing system and the proposed system including the number of stop points, number of trucks, time spent on the system, and distance traversed by trucks. The comparison shows that the proposed system would require less total travel time, but more total traveled distance. However, this is because of the increased number of trucks with less capacity than in the current system, with the stop points only at dumpster locations. While collecting waste door-to-door implies a more traveled distance and a higher number of stop points. However, the truck drivers in the current operating collection system do not implement the shortest path technique to find the optimal routes, thus the current actual travel time and travel distance in the field are believed to be significantly larger than the estimated numbers from the GIS model. In addition, the small trucks used in the proposed system are essential to allow for their accessibility to all stop points that would be required to fulfill the needs of the "door-to-door" collection system, and not contaminate the surroundings with large amounts of waste traveling long distances compared to the smaller trucks.

Therefore, the advantage of the proposed scenario lies in minimizing total travel time or the waste collection process to prevent waste dumping on the streets, which has a negative impact on the environment. The following section will quantify the environmental benefits of the proposed system.

## 4.2 Quantification of Environmental Impacts

Time and energy consumption during the waste collection procedure are the main components of the total cost of the waste management systems together with other environmental and economic impacts. Sonesson (2000) developed empirical equations that estimate both time and energy consumption. These equations generalize the calculations for waste management systems with a good outcome between a 5 and 24% deviation for time and energy consumption calculations. The developed equations calculate fuel consumption related to distance driving including traffic-related stops and extra stops due to collecting waste from stop points. Similarly, time consumption including hauling time, driving time, and emptying bins time is calculated.

Table 2 shows a comparison of the results between the existing and proposed scenarios, calculated using COPERT software. Both average route distance and time are decreased due to the use of smaller trucks in the proposed scenario. On the other hand, the number of trucks in the proposed scenario increased, and thus the total energy consumption is larger. We can justify the increase in total energy consumption by the positive environmental effect of the 'door-to-door' collection system compared to the collection from dumpsters. The 'door-to-door' collection system also ensures that the generated waste from households maintains its value before entering the treatment procedures, and thus generates more income, whilst in the other collection system it is prone to lose its value by scavengers who search for valuable items in the dumpsters, as well as MSW is mixed with street sweepings, which greatly affects the quality of MSW and suitability for treatment.

Table 2: Comparison between Existing and Proposed Scenario

| Scenario | Existing Scenario | Proposed Scenario | % Improvement |
|---|---|---|---|
| Waste Collection Method | Designated Dumpster | Door-to-Door | - |
| Number of Trucks | 16 | 50 | - |
| Truck Capacity (ton) | 18 | 4 | |
| Number of Stop points | 381 | 500 | - |
| Average Time spent at each Collection Point (min.) | 15 | 30 | - |
| Average Route Distance (km) | 110 | 67 | 39% |
| Total Traveled Distance (km) | 1,756 | 3,347 | - |
| Average Route Time (hr.) | 5.3 | 1.2 | 77% |
| Total Energy Consumption (MJ/day) | 108,907 | 336,089 | - |
| Total Time Consumption (h/day) | 84.6 | 62.2 | 26.5% |
| CO Emissions (g/day) | 142 | 97 | 32% |
| $CO_2$ Emissions (g/day) | 34,197 | 19,462 | 43% |
| $NO_x$ Emissions (g/day) | 473 | 210 | 56% |

## 5. Conclusions

In this study, the Network Analyst tool in ArcMap was used to plan an optimized and efficient waste collection and transportation plan on the 10[th] of Ramadan City in Cairo, Egypt. The plan is given a set of inputs including the collection fleet characteristics, proposed treatment, and disposal facility, proposed "door-to-door" collection system, and working time window. Both the planned system and the current operating system in the city are modeled in the ArcMap GIS platform and evaluated against each other. At first, the GIS City model is loaded into ArcMap and the Network Analyst tool is initialized through it. Then, the waste collection bins/trucks' stop point locations are set inside the city to achieve the criteria of service coverage of all residential areas. Also, the collection trucks are defined in the network, and the model is run to find the optimal routes and the minimum number of trucks that traverse all stop points in the city covering the collection bins/stop points within the defined working time window. Results show that the new proposed system improved the waste collection procedure by achieving the required task within less total travel time (26.5%) than the current system. In addition, the proposed collection system decreased the negative environmental impacts represented by CO, $CO_2$, and $NO_x$ emissions by (32%, 43%, and 56%) respectively. However, the total distance traveled in the proposed system is larger than the current system by 90% more, and that is due to the smaller size of used trucks (4 tons) than those used in the current system (18 tons) which increases the total traveling distance by more trucks.